\documentclass{article}

\usepackage{amsmath}
\usepackage{amsfonts}

\newcommand{\ms}{\medskip}

\newcommand{\noi}{\noindent}

\newcommand{\q}{{\cal Q}}
\newcommand{\bb}{{\cal B}}
\newcommand{\oo}{{\cal O}}
\newcommand{\fd}{finite-dimensional}
\newcommand{\id}{infinite-dimensional}
\newcommand{\fc}{finite-codimensional}
\newcommand{\pa}{Poisson algebra}

\newcommand{\la}{Lie algebra}

\newcommand{\ba}{basic algebra}
\newcommand{\ci}{C^{\infty}(M)}
\newcommand{\fb}{{\mathfrak b}}
\newcommand{\fa}{{\mathfrak a}}
\renewcommand{\ss}{{semisimple}}

\def\endproof{\hfill\vrule height4pt width6pt depth2pt}
\def\g{{\mathfrak g}}

\newtheorem{thm}{Theorem}
\newtheorem{cor}{Corollary}
\newtheorem{prop}{Proposition}
\newtheorem{defn}{Definition}

\title{An Obstruction to Quantizing Compact Symplectic Manifolds
\footnotetext{This is a revised version of the article 
published in {\sl Proc. Amer. Math. Soc.} {\bf 128}, 237--243 \\
\indent \indent \indent [2000]}}

\author{{\bf Mark J. Gotay}\thanks{Supported in part by NSF grant
96-23083. E-mail: gotay@math.hawaii.edu}
\\  Department of Mathematics \\ University of Hawai`i \\ 2565 The
Mall \\ Honolulu, HI 96822  USA \\ 
\and {\bf Janusz Grabowski}\thanks{Supported by KBN, grant No. 2 PO3A
042 10. E-mail: jagrab@mimuw.edu.pl} \\ Institute of
Mathematics
\\ University of Warsaw \\ ul. Banacha 2 \\ 02-097 Warsaw, Poland \\
\and {\bf Hendrik B. Grundling}\thanks{Supported by a research grant
from the Australian Research Council. E-mail:
hendrik@ \newline
\indent \indent \indent maths.unsw.edu.au} \\ Department of Pure
Mathematics
\\ University of New South Wales \\ P.O. Box 1 \\ Kensington, NSW
2033 Australia}

\date{May 30, 1997; Revised October 14, 1999}

\begin{document}


\maketitle

\begin{abstract} We prove that there are no nontrivial \fd\ Lie
representations of certain Poisson algebras of polynomials on a compact
symplectic manifold. This result is used to establish the
existence of a universal obstruction to quantizing a compact
symplectic manifold, regardless of the dimensionality of the
representation.

\end{abstract}


\begin{section}{Introduction}

\addtocounter{thm}{-1}
\addtocounter{cor}{-1}

Let $(M,\omega)$ be a symplectic manifold. Consider the
Poisson algebra
$C^\infty_c(M)$ of all compactly supported smooth functions on $M.$ In
a recent paper \cite{GM}, Ginzburg and Montgomery observed:
\begin{thm} 
\label{thm:gm} There exists no nontrivial finite-dimensional Lie
representation of
$C^{\infty}_c(M).$
\end{thm}

When $M$ is compact, this yields the ``no-go'' result:
\begin{cor} There exists no nontrivial finite-dimensional full
prequantization of $M.$
\label{cor:gmp}
\end{cor} 

\noi (The terminology is explained in the following section.) Although
not surprising on mathematical grounds, since
$C^\infty_c(M)$ is ``huge,'' this corollary does have physical
import, as one expects the quantization of a compact phase
space to yield a
\emph{finite-}dimensional Hilbert space.

{}From the standpoint of quantization theory, however, one is more
interested in certain polynomial algebras than in
$C^\infty_c(M)$ or even all of
$\ci.$ Inspired by this work, as well as that of Avez \cite{av1}, here
we generalize both Theorem
\ref{thm:gm} and Corollary \ref{cor:gmp} as follows. Let $\cal B$ be a
basic set of observables on $M$, $\fb$ the \la\ generated by $\bb$,
and
$P(\fb)$ the \pa\ of polynomials generated by $\fb$. (See \S2.) Assume
that $M$ is compact and that $\fb$ is finite-dimensional. Our main
result is:
\begin{thm}
\label{thm:gg} There exists no nontrivial
finite-dimensional Lie representation of
$P(\fb)$.
\end{thm}

\noi {}From this we immediately obtain:
\begin{cor}
\label{cor:pq} There exists no nontrivial
finite-dimensional prequantization of
$P(\fb)$.
\end{cor}

Furthermore, using Corollary \ref{cor:pq}, we are able to establish
the following result, which holds \emph{regardless} of the
dimensionality of the representation.
\begin{thm}
\label{thm:ggg} There exists no nontrivial
quantization of the pair $\left(P(\fb),\bb\right)$.
\end{thm}

Thus there is a universal obstruction to the quantization of a
compact symplectic manifold. Theorem \ref{thm:ggg} generalizes
some of the results obtained in \cite{GGH} for the phase space
$S^2$, with $\mathcal B = \mbox{su}(2)$, corresponding to a 
spinning particle. In this instance the polynomial algebra
$P(\mbox{su}(2))$ consists of spherical harmonics. Implications
of these no-go results will be discussed in the sequel.

\end{section}


\begin{section}{On Basic Algebras and Quantization}
 
Suppose that $(M,\omega)$ is a connected $2n$-dimensional symplectic
manifold with corresponding Poisson bracket \{\ ,\ \}. We first recall
from
\cite{GGT} the notion of a ``basic set of observables'' on $M$
which is central to our results. 

\begin{defn} {\rm \, A {\sl basic set of observables\/} $\bb$
is a \fd\ linear subspace of $\ci$ such that:\footnote{Unlike in
\cite{GGT} we do not require that $1 \in \bb$; this is superfluous.}
\begin{description} \item \rule{0pt}{0pt} \vspace{-18pt}
\begin{enumerate}
\item[(B1)] the Hamiltonian vector fields
$X_f,{f\in\bb}$, are complete, 
\vskip 6pt
\item[(B2)] $\bb$ is ``transitive,'' i.e. the collection
$\{X_f,{f\in\bb}\}$ spans $TM$, and
\vskip 6pt
\item [(B3)] $\bb$ is a minimal space satisfying these
requirements.
\end{enumerate} \end{description}}
\end{defn}

Given a basic set $\bb$, we denote by $\mathfrak b$ the Lie algebra
generated by $\bb.$ We call $\fb$ a {\sl basic algebra}. Throughout
this paper, we assume that
$\mathfrak b$ is \emph{finite}-dimensional. We refer the reader to
\cite[\S\S1.2.5-6]{on} for the necessary background on Lie algebras,
Lie groups, and their actions.

When $M$ is compact, we have the following characterization of basic
algebras. 
\begin{prop}
\label{prop:bs} 
Let $\fb$ be a \fd\ \ba\ on a compact symplectic manifold.
Then, as a \la, $\fb$ is compact and semisimple. In
particular, its center must be zero. 
\end{prop}

\noi \emph{Proof}. Define an inner product on $\fb$ according to
\begin{equation}
\langle f,g \rangle = \int_M fg\,\,\omega^n.
\label{eq:ip}
\end{equation}

\noi Using the identity
\begin{equation}
\label{eq:identity}
\{f,g\}\,\omega^n = n\, d(f\,dg \wedge \omega^{n-1})
\end{equation}

\noi together with Stokes' Theorem, we immediately verify that
\[\langle \{f,g\},h \rangle + \langle g,\{f,h\} \rangle = 0\]

\noi whence $\fb$ is compact. 
As a consequence,
$\fb$ splits as the direct sum $\mathfrak s \oplus \mathfrak z$, where
$\mathfrak s$ is semisimple and $\mathfrak z$ is the center of $\fb$.
Furthermore, the span of any generating set of
$\fb$ must then contain the center $\mathfrak z.$

Now the transitivity property (B2) implies that any function which
Poisson commutes with every element of $\fb$ must be a constant, so
that $\mathfrak z \subset \mathbb R.$ But if $\mathfrak z = \mathbb
R$ then, since $\mathfrak z \subset \bb$, $\bb$ would
not be minimal, thereby violating (B3). Thus
$\mathfrak z = \{0\}$ and $\fb$ is semisimple.
\endproof

\ms

Although we make no use of it, we remark that the converse of
Proposition~\ref{prop:bs} holds: Every \fd, compact,
semisimple \la\ $\fb$ is a basic algebra on some symplectic manifold
$M$. In fact, we can take $M$ to be a coadjoint orbit of maximal
dimension in $\fb^*$, cf. the discussion following the proof of
Theorem~\ref{thm:con} in \S 4.

\ms 

We denote by $P(\fb)$ the polynomial algebra generated by
$\fb$. Since $\fb$ is a Lie algebra, $P(\fb)$ is a Poisson algebra.
Note that (\emph{i})
$P(\fb)$ is not necessarily free as an associative algebra (cf. the
examples in
\S4), and (\emph{ii}) by definition $\mathbb R \subset P(\fb)$. Let
$P^k(\fb)$ denote the (\fd) subspace of polynomials of minimal
degree at most $k$. (Since $P(\fb)$ is not necessarily free, the
notion of ``degree'' is not well-defined, but that of ``minimal
degree'' is.)

\ms

Finally, we state what it means to quantize a Lie
subalgebra of observables. Again, we refer the reader to
\cite{GGT} for motivation and further discussion.

\begin{defn}$\,\,$ {\rm Let $\oo$ be a Lie subalgebra of
$C^\infty(M)$ containing $1$. A {\sl
prequantization\/} of $\oo$ is a Lie representation 
$\q$ of $\oo$ by skew-symmetric operators which preserve a common
dense domain $D$ in some separable Hilbert space such that 

\begin{description} \item \rule{0pt}{0pt} \vspace{-18pt}
\begin{enumerate}
\item[(Q1)] ${\cal Q}(1) = I$, and
\vskip 6pt
\item[(Q2)] if the Hamiltonian vector field of $f \in \oo$ is
complete, then
$\q(f)$ is essentially skew-adjoint on $D$.
\end{enumerate}
\end{description}

\noi If $\oo = C^\infty(M)$, the prequantization is said to be {\sl
full\/}.} 
\end{defn}

Now suppose that $\bb \subset
\oo$ is a basic set of observables.
\begin{defn} \, {\rm  A {\sl quantization\/} of the pair
$({\oo},{\bb})$ is a prequantization $\q$ of $\oo$ on a domain $D$
such that
\begin{description} \item \rule{0pt}{0pt} \vspace{-18pt}
\begin{enumerate}
\item[(Q3)] $\q(\bb) = \left\{\q(f)\,|\,f \in \bb\right\}$ is an
irreducible set, and 
\vskip 6pt
\item[(Q4)] $D$ contains a dense set of separately analytic vectors
for $\q(\bb).$
\end{enumerate}
\end{description}}
\label{defn:quant}
\end{defn}

When $\fb$ is \fd, (Q4) guarantees
that the representation of $\fb$ so obtained is
integrable to a representation of the corresponding simply
connected Lie group. We are interested here in the case when $\oo = P(\fb).$

\end{section}


\begin{section}{Proofs}

We now prove Theorem \ref{thm:gg}, Corollary
\ref{cor:pq}, and Theorem \ref{thm:ggg}. We begin with a purely
algebraic lemma. Let
$P$ be an abstract Poisson algebra, i.e. an associative commutative
algebra, equipped with a Poisson bracket
$\{\ ,\ \}$. The Poisson bracket is a Lie bracket satisfying the
Leibniz rule
\begin{equation*}
\{ f,gh\} =\{ f,g\} h+g\{ f,h\}.
\end{equation*} 

\noi The Leibniz rule implies the identity
\begin{equation}\label{2}
\{ f,gh \} +\{ g,hf\} +\{ h,fg\} =0.
\end{equation}
\begin{prop}
\label{prop:g} If $L$ is a finite-codimensional Lie ideal of a Poisson
algebra $P$ with identity, then either $L$ contains the derived ideal
$P'=\{ P,P\}$ or there is a maximal finite-codimensional 
associative ideal
$J$ of $P$ such that $P'\subset J$. \end{prop}

\noi {\em Proof.} Let us assume that $P'$ is not included in $L$ and
let us set $U_L=\{ f\in P\mid \{ f,P\}\subset L\}$.
\par {}From (\ref{2}) it follows easily that $U_L$ is an associative
subalgebra of $P$ (cf. \cite[1.6(b)]{Gr2}). Moreover,
$U_L$ is finite-codimensional, since it contains $L$, and
$U_L\ne P$, since the derived ideal $P'$ is not included in $L$.
\par Let $K$ be the largest associative ideal of $P$ included in
$U_L$, i.e. $K=\{ f\in P\mid fP\subset U_L\}$. As the kernel of the
regular representation of $U_L$ in $P/U_L$ by associative
multiplication, $K$ is finite-codimensional in $U_L$ (and hence in
$P$) and $K\ne P$. Moreover, $\{ P,K\}\subset L\subset U_L$. Since for
any $f\in P$ we have
$$ f\{ P,K\}\subset \{ P,fK\} +\{ P,f\} K\subset U_L,$$ we can write
$\{ P,K\}\subset K$.
\par Since $P$ has an identity, $K$ is contained in some maximal
(and thus prime) \fc\ associative ideal $J.$ Let
$X_f$ be the Hamiltonian vector field of $f\in P$. Since $$ J^{X_f}=\{
g\in P\mid (X_f)^ng\in J
\quad{\rm for\ all}\ n=0,1,2,\dots\}$$ is finite-codimensional (as it
contains $K$), we conclude from \cite[Lemma 4.1]{Gr1} that
$X_f(J)\subset J$, i.e.
$\{ P,J\}\subset J$. Finally, we use \cite[Lemma 4.2]{Gr1} to get that
$\{ P,P\}\subset J$.
\endproof

\ms

Since the basic algebra
$\fb$ is fixed, we will henceforth abbreviate the \pa\ $P(\fb)$ by
$P$. 

\ms

\emph{Proof of Theorem \ref{thm:gg}.} Suppose that $\q$ were a Lie
representation of $P$ on some \fd\ vector space. Then $L = \ker \q$
is a \fc\ Lie ideal of $P.$ We will show that $L$ has codimension at
most 1, whence the representation is trivial (i.e. $\q$  factors
through a representation of $P/L$ with $\dim P/L \leq 1.$) We
accomplish this in two steps, by showing that:
\begin{description} \item \rule{0pt}{0pt} \vspace{-18pt}
\begin{enumerate}
\item[(a)] The derived ideal $\{P,P\}$ has codimension 1 in
$P,$ and
\item[(b)]  $L \supset \{P,P\}$.
\end{enumerate}
\end{description}

Recall that the mean of
$f \in \ci$ is
\[\bar f = \frac{1}{\mbox{vol}(M)}\,
\int_M f\,\omega^n.\] Let $P_0$ denote the Lie
ideal of all polynomials of zero mean. The decomposition $f
\mapsto \bar f + (f -
\bar f)$ gives 
$P = \mathbb R \oplus P_0$. Thus, if we prove that $\{P,P\} = P_0$,
(a) will follow.

Using (\ref{eq:identity}) along with Stokes' Theorem, we immediately
have that
$\{P,P\}
\subset
P_0$. To show the reverse inclusion, let
$\{b_1,\ldots,b_N\}$ be a basis for
$\fb$, so that
\[\{b_i,b_j\} = \sum_{k=1}^N c^k_{ij}b_k\]

\noi for some constants $c^k_{ij}$. Following Avez \cite{av2},
define the ``symplectic Laplacian''
\[\mathnormal{\Delta}f = - \sum_{i = 1}^N \,\{b_i,\{b_i,f\}\}.\]

\noi It is clear from these two expressions and the Leibniz rule that
the linear operator
$\mathnormal{\Delta}$ maps
$P^k$ into $P_0^{\,k}.$ Furthermore, taking into account the
transitivity assumption (B2), we can apply
\cite[Prop. 1(4)]{av2} to conclude that
$\mathnormal{\Delta}f = 0$ only if $f$ is constant. Thus for each $k
\geq 0$, the decomposition $P^k = \mathbb R \oplus P_0^{\,k}$ implies
$\mathnormal{\Delta}(P^k) = P_0^{\,k}.$ It follows that $P_0 \subset
\{P,P\}.$

\ms 

If (b) does not hold, then by Proposition \ref{prop:g} there must be a
proper associative ideal $J$ in $P$ with $\{P,P\} \subset J.$
Since
$\{P,P\}= P_0$ has codimension 1, $P_0 = J.$ This is, however,
impossible, since $f^2$ has zero mean only if $f=0. $ \endproof

\ms

Corollary \ref{cor:pq} is trivial, since a prequantization is simply a
special type of Lie representation.

\ms 

\emph{Proof of Theorem \ref{thm:ggg}}. We proceed by
reducing to the \fd\ case and applying Corollary \ref{cor:pq}.
Suppose that $\q$ were a quantization of
$(P,\bb)$ on a Hilbert space. By conditions (Q2) and (Q4), $\q(\fb)$
can be exponentiated to a unitary representation of the connected,
simply connected Lie group $B$ with \la\ $\fb$ (recall that $\fb$ is
assumed \fd) which, according to (Q3), is irreducible. Since by
Proposition \ref{prop:bs} $\fb$ is compact and semisimple, $B$ is
compact.
The representation space must thus be
\fd, and so Corollary \ref{cor:pq} applies.
\endproof

\end{section}


\begin{section}{Discussion}

Theorem \ref{thm:ggg} asserts that the polynomial algebra $P(\fb)$
generated by any \fd\ basic algebra $\fb$ on a compact
symplectic manifold cannot be consistently quantized. We emphasize,
however, that this result need not hold if $\fb$ is allowed to be
in\fd: there exists a full quantization of the
torus \cite{Go}. (In this regard we observe that by
Theorem \ref{thm:char} below, there is no \emph{finite-}dimensional
\ba\ on the torus as it is not simply connected.) Similarly, Theorem
\ref{thm:gg} and Corollary
\ref{cor:pq} can fail when the  representation space is allowed to be
\id: as is well-known, full prequantizations exist provided
$\omega$ is integral. Thus Corollary \ref{cor:pq} and Theorem
\ref{thm:ggg} are the optimal no-go results for compact phase spaces.

It is of interest to determine those compact symplectic manifolds
which admit basic sets of observables (or, equivalently, \ba s). It
turns out that these form a quite restricted class. 
\begin{thm} If a compact, connected, symplectic manifold
$M$ admits a finite-di\-men\-sion\-al
\ba\
$\fb$, then $M$ is a coadjoint orbit in $\fb^*.$
In particular, $M$ must be simply connected.
\label{thm:char}
\end{thm}
\emph{Proof.} Let $\fb$ be a \ba\ on $M,$ and denote by $\g$ the
\la\ formed by the Hamiltonian vector fields of elements of
$\fb.$ Since $M$ is compact the vector fields in $\g$ are complete, so
$\g$ may be integrated to a symplectic action on $M$ of some
connected Lie group $G$ with Lie algebra $\g$.
Condition (B2) implies that this action is
locally transitive, and therefore globally transitive as $M$ is
connected. Thus $M$ is a homogeneous space for $G$.

By
Proposition \ref{prop:bs} the center of $\fb$ is trivial. It follows
that $\g
\cong \fb$, whence the action of $G$ on $M$ is Hamiltonian. Again
using the compactness and semisimplicity of
$\fb$, we conclude 
that $G$ is compact. 
The Kirillov-Kostant-Souriau Coadjoint Orbit Covering Theorem
\cite[Thm. 14.6.5]{m-r} then implies that $M$ is a symplectic covering
of a coadjoint orbit $O$ in $\fb^*$. On the other hand, the
coadjoint orbits of a compact connected Lie group are simply
connected \cite[Thm. 2.3.7]{Fi}, so that $M$ is symplectomorphic
to $O$. 
\endproof 

\ms
As a partial converse to Theorem~\ref{thm:char}, we have
\begin{thm} Let $M$ be a nonzero coadjoint orbit in $\fb^*$, where
$\fb$ is a compact, simple \la. Then $M$ admits $\fb$ as a \ba.
\label{thm:con}
\end{thm}
\emph{Proof.} We first observe that since $\fb$ is compact and
semisimple, its adjoint group $B$ is compact \cite[Thm.
4.11.7]{Va}.
Consequently $M$ is compact. 

Now the elements of $\fb$, regarded as
(linear) functions on
$\fb^*$, form a finite-dimensional space of observables on $M$ closed
with respect to the Poisson bracket $\{\ ,\ \}$.
Their Hamiltonian vector fields are complete due to the
compactness of $M$; moreover, transitivity is automatic since $M$ is
an orbit. Thus $\fb$ satisfies (B1) and (B2). We will
show below that
$\fb$ is minimal with respect to these properties. Given this, let
$\bb$ be a subspace of $\fb$ which is minimal amongst all
such satisfying (B2). Then $\bb$ is a basic set and generates $\fb$, 
for otherwise the transitive \la\ generated by $\bb$ would be smaller
than
$\fb$, which contradicts the minimality of $\fb$. It follows that
$\fb$ is a basic algebra.

Identifying $\fb$ with $\fb^*$ by means of the {\em Ad\/}-invariant
inner product (1), we can identify adjoint
and coadjoint orbits. Thus $M = O_h$ for some $h \in \fb.$ With
respect to this identification
$\{ f,g\}(k)= \langle [f,g],k \rangle$ for
$f,g,k\in \fb$. It is easy to see that the isotropy subalgebra
$\fb_h=\{ f\in
\fb\,|\,[f,h]=0\}=[\fb,h]^\perp$, so we have
the decomposition 
\setcounter{equation}{3}
\begin{equation}
\fb=\fb_h\oplus [\fb,h]
\label{decomp}
\end{equation}

\noi with
\begin{equation}
[\fb_h,[\fb,h]]\subset [\fb,h].
\label{eq:g0}
\end{equation}

Suppose now that $\fb$ is not minimal, and let $\fa \subset \fb$ be
a minimal subalgebra for the compact orbit $O_h$.  Proposition
1 implies that $\fa$ is compact and \ss , and hence has a
commutative Cartan subalgebra
$\mathfrak c$. We extend $\mathfrak c$ to a Cartan subalgebra
$\mathfrak k$ of $\fb$. Now fix a Cartan subalgebra of
$\fb$ containing $h$; since this subalgebra is commutative it must
be contained in $\fb_h.$ As all Cartan
subalgebras of a compact \ss\ \la\ are conjugate under $B$ \cite[Thm.
4.12.2]{Va}, we may map $\mathfrak k $ onto this subalgebra by some
conjugation. Since
$O_h$ is an adjoint orbit, the image of $\fa$ under this conjugation
is also basic for
$O_h$. Thus without loss of generality we may suppose that
$\mathfrak c \subset \mathfrak k \subset \fb_h$.

Passing to complexifications we have the root space decomposition
$$\fa_\mathbb C = \mathfrak c_\mathbb C \oplus \left (\; \sum_{\delta \in
\Delta}\fa_\delta\right )$$
where $\Delta \subset \mathfrak c_\mathbb C^{\;\;*}$ is the set of roots
of
$(\fa_\mathbb C,\mathfrak c_\mathbb C).$ The inclusion 
$\mathfrak c \subset \mathfrak k$ 
implies that every root space $\fa_\delta$ of
$(\fa_\mathbb C,\mathfrak c_\mathbb C)$ is also a root space
$\fb_{\delta '}$ of $(\fb_\mathbb C,\mathfrak k_\mathbb C)$ for some unique
$\delta' \in
\mathfrak k_\mathbb C^{\;\; *}$ with 
$\delta' \, |\, \fa = \delta$. 
We can therefore write
\begin{equation}
\fa_\mathbb C = \mathfrak c_\mathbb C \oplus \left (\; \sum_{\delta
\in \Delta}\fb_{\delta '}\right ).
\label{rs}
\end{equation}

\noi By transitivity $[\fa,h] = T_hO_h = [\fb,h]$. Recalling that $h
\in \mathfrak k$ so that $[\fb_{\delta'},h] \subset \fb_{\delta '}$, we
compute using (\ref{rs}),
\begin{equation*}
[\fb_\mathbb C,h] = [\fa_\mathbb C,h] = \left[\mathfrak c_\mathbb C,h \right ]
\oplus
\left (\; \sum_{\delta \in \Delta}\left [\fb_{\delta '},h  \right ]
\right )  \subset 
\sum_{\delta \in \Delta}\fb_{\delta '}
  = \sum_{\delta \in \Delta}\fa_{\delta}\subset \fa_\mathbb C.
\end{equation*}

\noi Thus 
\begin{equation}
[\fb,h] \subset \fa.
\label{=}
\end{equation}

Let $\mathfrak h$ be the \la\ generated by $[\fb,h]$; then from
(\ref{decomp}) and (\ref{eq:g0}) we see that
$\mathfrak h$ is an ideal in $\fb$. From (\ref{=}) it follows that
$\mathfrak h \subseteq \fa.$ But as $\fb$ is simple, $\mathfrak h =
\fb$ and this forces $\fa = \fb$, whence $\fb$ is minimal.
\endproof

\ms

When $\fb$ is merely semisimple this
result need not hold. For instance, $S^2$ is a coadjoint orbit in
$\mbox{su(2)}^* \oplus \{0\} \subset \mbox{su(2)}^* \oplus
\mbox{su(2)}^*
\cong \mbox{so(4)}^*$ with basic algebra $\mbox{su(2)}$, not
$\mbox{so(4)}.$ (More generally, if there are several simple
components, say
$\fb=\fb_1\oplus \cdots \oplus \fb_K$, then there will be
low-dimensional orbits for which $\fb$ is not
basic; for instance, for
$h$ a regular element in a Cartan subalgebra $\mathfrak c$ of $\fb_1$
we have
$\fb_h=\mathfrak c\oplus
\fb_2 \oplus \cdots \oplus \fb_K$, and the corresponding orbit is
$B/B_h \cong B_1/C$ with basic algebra
$\fb_1$, not
$\fb$.)

When $M$ has \emph{maximal} dimension in $\fb^*$, however, $\fb$ will
be a basic algebra on $M$. Indeed, 
it is a well-known result of Duflo and Vergne that the isotropy
algebras of maximal coadjoint orbits are abelian (see, e.g.
[{\bf MR}, pps. 278-280]). But, referring back to the proof of
Theorem~\ref{thm:con} it follows from this fact, (\ref{decomp}),
(\ref{eq:g0}), and the semisimplicity of
$\fb$ that $[\fb,h]$ generates $\fb$, i.e. $\mathfrak h = \fb$.
So again (\ref{=}) yields $\fa=\fb$.


\ms

Given that 
$M$ is a coadjoint orbit in
$\fb^*$, we may identify
$P(\fb)$ with the Poisson algebra of polynomials on
$\fb^*$ restricted to $M$. In particular, we can take $M=S^2 \subset
\mbox{su(2)}^*$,
$\bb$ the space of spherical harmonics of degree one ($\bb = \fb \cong
\mbox{su}(2)$), and
$P(\fb)$ the space of all spherical harmonics. In \cite{GGH} it was
shown with some effort that there is no nontrivial quantization of
the pair
$\left(P(\fb),\bb\right).$ But this now follows immediately from
Theorem \ref{thm:ggg}. A similar analysis applies to $\mathbb CP^n
\subset
\mbox{su}(n+1)^*$, thereby answering in the affirmative a question
posed in
\cite[\S3]{GGH}.

Finally, we remark that the map $f \mapsto \bar f$ always
provides a trivial--but nonzero--prequantization of $\ci$ on
$\mathbb C$. When the representation space is in\fd, Avez \cite{av1}
has shown that $f \mapsto \bar fI$ is the only
possible prequantization of $\ci$ by
\emph{bounded} operators.

We explore the case when $M$ is noncompact in a separate work
\cite{GG2}.

\end{section}


\section*{Acknowledgments}

We would like to thank V. Ginzburg and G. Folland for helpful
discussions.



\end{document}